# ADAPTIVE NONPARAMETRIC CONFIDENCE SETS

By James Robins and Aad van der Vaart

*Harvard University and Vrije Universiteit Amsterdam*

We construct honest confidence regions for a Hilbert space-valued parameter in various statistical models. The confidence sets can be centered at arbitrary adaptive estimators, and have diameter which adapts optimally to a given selection of models. The latter adaptation is necessarily limited in scope. We review the notion of adaptive confidence regions, and relate the optimal rates of the diameter of adaptive confidence regions to the minimax rates for testing and estimation. Applications include the finite normal mean model, the white noise model, density estimation and regression with random design.

**1. Introduction.** Consider an observation $X^{(n)}$ distributed according to a law $P_\theta^{(n)}$ depending on a parameter $\theta$ that ranges over a subset $\Theta$ of a separable Hilbert space. Specifically, we take the Hilbert space equal to $\mathbb{R}^n$ with the Euclidean norm, or the sequence space $\ell_2 = \{\theta = (\theta_1, \theta_2, \dots) : \sum_{i=1}^\infty \theta_i^2 < \infty\}$ with the squared norm $\|\theta\|^2 = \sum_{i=1}^\infty \theta_i^2$. Our aim is to construct (asymptotic) confidence sets $\hat{C}_n$ of small diameter for the parameter $\theta$, which are "honest" in the sense that, for a given confidence level $1 - \alpha$,

$$(1.1) \qquad \liminf_{n \to \infty} \inf_{\theta \in \Theta} P_\theta(\theta \in \hat{C}_n) \geq 1 - \alpha.$$

This problem has been considered by, among others, Li [32] and Baraud [1] in the case that $\Theta$ is equal to $\mathbb{R}^n$ and the observation is a Gaussian vector with mean $\theta$ and covariance matrix the identity, by Hoffmann and Lepski [20] in the case that $\theta \in \ell_2$ and the observation is an infinite sequence of Gaussian variables with means $\theta_i$ and variance $\sigma^2/n$, and by Beran [4], Beran and Dümbgen [5] and Genovese and Wasserman [18] in the case of the fixed design regression model. Our aim in this paper is to propose new confidence procedures for these and related models, which shed light on some of the questions raised in the discussion of the paper by Hoffmann and Lepski [20]. We construct confidence sets with the properties:

(i) The confidence set is honest on the model $\Theta$.

(ii) The confidence set is centered at an estimator of choice, for example, an adaptive estimator.

(iii) The diameter of the confidence set adapts to submodels of $\Theta$ in a rate-optimal way.

In the second and third points we improve on the results in the mentioned papers, at least as regards rates. Our method in its simplest form as presented below leads to an increase of the "constants."

Since completing our paper we have learned about the work of Juditsky and Lambert-Lacroix [25] and Cai and Low [11]. Juditsky and Lambert-Lacroix [25] appear to deserve priority in discussing adaptive confidence sets. In their beautiful paper they pose the problem within the setting of fixed-design regression with Gaussian errors and obtain adaptation in the scale of Besov spaces, using wavelet-based methods. An insightful discussion of the problem and basic insights about its relationship to loss estimation and minimax estimation and testing can already be found in this paper. Cai and Low [11] consider the problem of adaptive confidence regions in the setting of the Gaussian white noise model, and obtain adaptation in the scale of Besov spaces, also using wavelet-based estimators. Our method is more flexible and applies to more settings, but we develop the results only for the scale of Sobolev spaces. In certain respects it is close to the method of Juditsky and Lambert-Lacroix [25].

As is pointed out in the preceding references, the desired honesty (i) severely limits the possibility of adaptation as in aim (iii). In the past years many successes have been obtained in the construction of estimators that are simultaneously minimax over a large selection of models. (See, e.g., [2, 3, 13, 14, 15, 16, 19, 29, 30, 31, 34, 37].) These estimators are able to adapt to the "regularity" of the true underlying parameter, without pre-knowledge of the parameter or its regularity. However, as pointed out by Birgé [7], these estimators have the property of being close to the true parameter without the statistician being able to tell how close it is. An adaptive estimator can adapt to an underlying model, but does not reveal which model it adapts to, with the consequence that nonparametric confidence sets are necessarily much wider than the actual discrepancy between an adaptive estimator and the true parameter.

If one drops "honesty" (i) from the requirements of the confidence set, but requires, for instance, only that the confidence set is honest over every submodel $\Theta_1 \subset \Theta$ of interest [i.e., (1.1) with $\Theta$ replaced by $\Theta_1$], then this embarrassing problem disappears, and it is possible to construct "confidence sets" of a diameter that adapts to the estimation rate. Most procedures in the literature fall in this category. However, dropping full honesty (i) appears to contradict the very definition of a confidence set. In this paper we require honesty in the sense of (1.1) with $\Theta$ the collection of all parameters deemed possible. Thus we consider a list of models and require honesty on the "biggest model" $\Theta$ in the list.

Under this requirement the possibilities for adaptation are severely limited. For a given submodel $\Theta_1 \subset \Theta$, the diameter of a confidence region that is honest for $\Theta$ cannot be of smaller order, uniformly over $\Theta_1$, than:

(a) The "slowest rate" $\varepsilon_n \to 0$ such that for any estimator sequence $T_n$ and some $\beta > \alpha$

$$\liminf_{n \to \infty} \sup_{\theta \in \Theta_1} \mathrm{P}_\theta(\|T_n - \theta\| \geq \varepsilon_n) > \beta.$$

This is typically the minimax rate of estimation for the model $\Theta_1$.

(b) The minimax rate of testing of the hypothesis $H_0 : \theta \in \Theta_1'$ versus the alternative $H_1 : \theta \in \Theta, \|\theta - \Theta_1'\| > \varepsilon_n$, for any given $\Theta_1' \subset \Theta_1$, for example, a one-point set $\Theta_1' = \{\theta_1\}$. This rate is often determined by the full model $\Theta$, rather than the submodel $\Theta_1$.

These lower bounds appear to be well known. Juditsky and Lambert-Lacroix [25] discuss such bounds in the setting of Besov spaces. For completeness we give precise statements in Section 6.

Our confidence sets have diameter of the order the maximum of the rates in (a)–(b), simultaneously over many submodels, at least for regularity classes as in the following example, and hence satisfy aim (iii).

EXAMPLE 1.1 (Regular parameters).   A parameter $\theta \in \ell_2$ can be called $\beta$-regular (for a given $\beta > 0$) if it belongs, for some $L > 0$, to the ellipsoid

$$S(\beta, L) = \left\{ \theta \in \ell_2 : \sum_{i=1}^{\infty} \theta_i^2 i^{2\beta} \leq L \right\}.$$

If the coordinates of $\theta$ correspond to classical Fourier coefficients, then $S(\beta, L)$ corresponds to periodic functions with $\beta$ derivatives bounded by a multiple of $L$ in $L_2[0,1]$. (For real functions and the sine–cosine basis the correspondence is more accurate if we replace $i^{2\beta}$ by $(i-1)^{2\beta}$ for odd values of $i$. See, e.g., the Appendix of [38].)

Consider inference on $\theta \in S(\beta, L)$ based on observing each $\theta_i$ with an independent $N(0, \sigma^2/n)$ error, or in one of the other models discussed below, which yield similar results. The minimax estimation rate for $S(\beta, L)$ is $n^{-\beta/(2\beta+1)}$ (cf., e.g., [9, 21, 22, 36]). For $\beta_1 > \beta$ and $L_1 \leq L$ we have $S(\beta_1, L_1) \subset S(\beta, L)$ and the minimax testing rate of $S(\beta_1, L_1)$ relative to $S(\beta, L)$ in the sense of (b) is $n^{-\beta/(2\beta+1/2)} \ll n^{-\beta/(2\beta+1)}$. (See [23], Theorem 2.1 or 3.1, or [24].)

Table 1
*Order of maximal diameter of confidence regions on the submodel $S(\beta_1, L_1) \subset S(\beta, L)$, cut-off points and number of observations needed to estimate $\sigma^2$*

| $\beta$ | $\beta_1$ | Radius on $S(\beta_1, L_1)$ | Cut-off | Obs for $\sigma^2$ |
|---|---|---|---|---|
| 1 | $\geq 2$ | $n^{-2/5}$ | $n^{2/5}$ | $\gg n^{2/5}$ |
| 1 | $3/2$ | $n^{-3/8}$ | $n^{2/5}$ | $\gg n^{2/5}$ |
| 1 | $1$ | $n^{-1/3}$ | $n^{2/5}$ | $\gg n^{2/5}$ |
| $1/2$ | $\geq 1$ | $n^{-1/3}$ | $n^{2/3}$ | $\gg n^{2/3}$ |
| $1/2$ | $3/4$ | $n^{-3/10}$ | $n^{2/3}$ | $\gg n^{2/3}$ |
| $1/2$ | $1/2$ | $n^{-1/4}$ | $n^{2/3}$ | $\gg n^{2/3}$ |
| $1/4$ | $\geq 1/2$ | $n^{-1/4}$ | $n$ | $\gg n$ |
| $1/4$ | $1/4$ | $n^{-1/6}$ | $n$ | $\gg n$ |
| $1/8$ | $\geq 1/4$ | $n^{-1/6}$ | $n^{4/3}$ | $\gg n^{4/3}$ |
| $0$ | $\geq 0$ | $1$ | $n^2$ | $\gg n^2$ |

If the supermodel $\Theta$ is equal to $S(\beta, L)$, then these bounds suggest that the diameter of a confidence set can be of diameter of order $n^{-\beta/(2\beta+1)}$ uniformly over $\Theta$, and of order $n^{-\beta_1/(2\beta_1+1)} \vee n^{-\beta/(2\beta+1/2)}$ uniformly over the smaller model $\Theta_1 = S(\beta_1, L_1)$ for $\beta_1 > \beta$ and $L_1 \leq L$. If $\beta_1 \in (\beta, 2\beta)$, then the latter rate is equal to $n^{-\beta_1/(2\beta_1+1)}$ and depends on the submodel. In that case we may say that adaptation occurs.

This type of adaptation is very different from adaptation in the context of estimation. For $\beta_1 \geq 2\beta$ the diameter is $n^{-\beta/(2\beta+1/2)}$, independent of the exact value of $\beta_1$, so that further regularity does not yield smaller confidence regions. Even on very small submodels ($\beta_1 \to \infty$), the diameter of a confidence region is at least of the order $n^{-\beta/(2\beta+1/2)}$, determined by the supermodel. As illustration, Table 1 gives the rates for some values of the regularity parameters. The meaning of the last two columns of the table is explained later on in the paper.

Our method to construct confidence regions, described in Section 2, is based on a sample-splitting procedure. We use half the data to construct centering estimators $\hat{\theta}^{(n)}$, and an independent second half to construct a confidence region around $\hat{\theta}^{(n)}$. The nature of the initial estimator $\hat{\theta}^{(n)}$ is irrelevant for the honesty of the confidence procedure, and hence $\hat{\theta}^{(n)}$ can be any of our favorite estimators. In particular, it can be an estimator that adapts to a selection of models of our choice. Our procedure borrows its adaptive strength from these initial estimators, but of course only up to the limitations described earlier.

Refinements of this procedure would be to construct two confidence sets, with the roles of the two half-samples interchanged, and to take the intersection, or to split the sample into more parts. For restricted supermodels $\Theta$ the splitting may be avoided altogether. This may lead to better constants in the centering and diameter of the confidence set. In this paper we are interested in rates only, and for this our simple sample-splitting scheme suffices.

In the case that the observations are a random sample, we can form the two halves of the data by simply splitting the sample into two parts, using the first half-sample to construct the estimator $\hat{\theta}^{(n)}$ and the second to construct the confidence region. In other examples of interest a similar situation can be created using a more involved splitting device, which we describe below.

The organization of the paper is as follows. In Section 2 we describe the construction in a general framework. In Sections 3, 4 and 5 we give the details for the three main examples, sequence models, density estimation and random regression. Finally in Section 6 we relate the diameter of a confidence region to the testing and estimation rates.

We close this introduction with a description of a number of examples to which our construction applies, together with a review of the literature.

EXAMPLE 1.2 (Finite sequence model). In this model the observation is a vector $X^{(n)} = (X_1, X_2, \ldots, X_n)$ from an $n$-dimensional normal distribution with mean vector $\theta = (\theta_1, \theta_2, \ldots, \theta_n)$ and covariance matrix $(\sigma^2/n)I$. The variance $\sigma^2$ is known and the parameter $\theta$ is known to belong to a subset $\Theta$ of $\mathbb{R}^n$, which may be all of $\mathbb{R}^n$.

This model was studied in [32] and [1] under the assumption that $\Theta = \mathbb{R}^n$. The naive procedure in this situation is the chi-square region $\{\theta \in \mathbb{R}^n : \|\theta - X^{(n)}\|^2 \leq (\sigma^2/n)\chi^2_{n,1-\alpha}\}$, which derives from inverting the likelihood ratio test. It has diameter of order 1, uniformly in (and independently of) $\theta$.

Li [32] showed that requiring honesty relative to all parameters $\theta \in \mathbb{R}^n$ implies that no confidence region can achieve a diameter that is uniformly smaller than $n^{-1/4}$, and exhibits confidence regions around shrinkage estimators that may achieve the rate $n^{-1/4}$ on the submodel where the shrinkage estimator performs well. Li's confidence sets improve on the naive chi-square procedure at true parameters where the shrinkage estimator improves upon the naive estimator $X^{(n)}$. Baraud [1] constructs confidence regions that improve on the naive procedure in a wider range of submodels. His procedure is based on comparing a range of submodels by chi-square tests. The confidence regions in the present paper manage to adapt to still more submodels, if the initial estimators $\hat{\theta}^{(n)}$ are chosen so as to fully profit from the recent insights in adaptive estimation, such as in [8].

It is notable that in this model the variance $\sigma^2$ is assumed known. Baraud [1] shows that in the case that $\sigma^2$ is an unknown parameter ranging over some interval (even a very short one), confidence regions that are honest over $\Theta = \mathbb{R}^n$ and $\sigma^2$ can never have diameter less than order 1.

Because the observations in this example are non-i.i.d., splitting the sample is not a good device in order to separate constructing a center and a radius of the confidence region. However, we may artificially produce two normal vectors $X'$ and $X''$ with means $\theta$ from a given $N_n(\theta, (\sigma^2/n)I)$-distributed random vector $X$ using randomization. Given a sample of independent, uniform variables $U_i$ independent of $X$, it suffices to define

$$X_i' = X_i + \Phi^{-1}(U_i)\sigma/\sqrt{n},$$
$$X_i'' = X_i - \Phi^{-1}(U_i)\sigma/\sqrt{n}.$$

Then it can be verified that $X_i'$ and $X_i''$ are independent random variables with means $\theta_i$ and variances $2\sigma^2/n$. Thus the observations can be duplicated at the cost of multiplying the variance $\sigma^2$ by 2. In the remainder of the paper we shall assume that a device of this type has been applied, and write $X^{(n)}$ for the second sample (on which the estimate of the radius of the confidence set is based), and assume that this is independent of the initial estimator $\hat{\theta}^{(n)}$ for $\theta$.

Knowledge of $\sigma^2$ is crucial for this randomization step. Good estimators would do as well, but it is impossible to estimate $\sigma^2$ in this model without restricting the mean parameter $\theta$ to a proper subset $\Theta$ of $\mathbb{R}^n$. Baraud [1] shows that the size of a confidence set can never be of smaller order than the imprecision in $\sigma$.

EXAMPLE 1.3 (Infinite sequence model). In this model the observations are an infinite sequence $X^{(n)} = (X_1, X_2, \dots)$ of independent random variables $X_i$ possessing normal distributions with means $\mathrm{E}X_i = \theta_i$ and variance $\sigma^2/n$. The parameter is the mean vector $\theta = (\theta_1, \theta_2, \dots)$ and is known to belong to a subset $\Theta$ of $\ell_2$.

This model is a version of the white noise model, and is considered in connection to confidence regions in Hoffmann and Lepski [20]. (The focus of these authors is on "random normalizing constants" rather than confidence regions, but, as most of the discussants of their paper, we interpret their results with respect to their implications for confidence regions.) Hoffmann and Lepski [20] assume that there is a largest model $\Theta$ of interest, and exhibit confidence regions that are adaptive to finitely many submodels. Our construction allows infinitely many submodels and yields confidence regions around arbitrary initial estimators $\hat{\theta}^{(n)}$, for example, adaptive ones. Hoffmann and Lepski consider the general setting of anisotropic regression models, but we illustrate our method for the regularity classes of Example 1.1 only.

We can use the same device as in Example 1.2 to duplicate the observations, at the cost of doubling the variance $\sigma^2$.

Typically one chooses $\Theta$ to be a relatively small subset of $\ell_2$. Then it is easy to find good estimators of $\sigma^2$, and it is not necessary to assume that $\sigma^2$ is a priori known. For instance, if $\Theta$ is an ellipsoid of the form $\{\theta \in \ell_2 : \sum_{i=1}^{\infty} \theta_i^2 i^{2\beta} < \infty\}$, then we may base an estimate of $\sigma^2$ on the observations $X_{k+1}, X_{k+2}, \ldots, X_{k+m}$ for sufficiently large integers $k, m$, which are approximately $N(0, \sigma^2)$-distributed for large $k$. The availability of an infinite sequence allows one to control the bias and variance of estimators of $\sigma^2$ to arbitrary precision by choosing $k$ and $m$, respectively, sufficiently large.

EXAMPLE 1.4 (Density estimation). In this model the observation is an i.i.d. sample $X_1, \ldots, X_n$ from a density $f$ relative to some measure $\mu$ on a measurable space $(\mathcal{X}, \mathcal{A})$. The density $f$ is known to belong to a subset $\mathcal{F}$ of $L_2(\mathcal{X}, \mathcal{A}, \mu)$.

We can cast this example into a problem of estimating a sequence $\theta = (\theta_1, \theta_2, \ldots)$ of parameters by expanding $f$ on a fixed orthonormal basis $e_1, e_2, \ldots$ of $L_2(\mathcal{X}, \mathcal{A}, \mu)$. This expansion takes the form of the Fourier series $f = \sum_i \theta_i e_i$, for the Fourier coefficients $\theta_i = \langle f, e_i \rangle = \mathrm{E} e_i(X_1)$.

The empirical Fourier coefficients $Y_i = n^{-1} \sum_{j=1}^{n} e_i(X_j)$ are unbiased estimators of the parameters $\theta_i$. However, they are only approximately normally distributed and not independent, and it seems not fruitful to cast this example into the framework of the sequence model of Example 1.3 with observational vector $(Y_1, Y_2, \ldots)$. The Le Cam equivalence of the white noise model and the density estimation model, proved under conditions by Nussbaum [35], offers a different connection between the two examples, but can be used only if $\mathcal{F}$ is restricted and yields regions of complicated form. (The latter objection is alleviated by the recent constructions of Brown, Carter, Low and Zhang [10].) Our direct approach gives concrete confidence sets and in wider generality.

We can split the sample into two independent halves to construct the center $\hat{\theta}^{(n)}$ and the radius $R_n(\hat{\theta}^{(n)})$ of the confidence set.

There is no parameter $\sigma^2$ to be dealt with in this example.

EXAMPLE 1.5 (Random regression). In this model the observation is an i.i.d. sample $(X_1, Y_1), \ldots, (X_n, Y_n)$ from the distribution of a vector $(X, Y)$ described structurally as $Y = f(X) + \varepsilon$, for $(X, \varepsilon)$ a random vector with $\mathrm{E}(\varepsilon | X) = 0$ and $\mathrm{E}(\varepsilon^2 | X) < \infty$ almost surely. The regression function $f$ is known to belong to a subset $\mathcal{F}$ of $L_2(\mathcal{X}, \mathcal{A}, P_X)$ for $P_X$ the marginal distribution of $X$, which is assumed known. The variance function $\sigma^2(x) = \mathrm{E}(\varepsilon^2 | X = x)$ need not be known, although for confidence intervals that are honest in $\sigma^2$ we need a known upper bound. We do not assume that the errors are normally distributed, and we do not assume that $X$ and $\varepsilon$ are independent.

As in Example 1.4 we can cast this example into a problem of estimating a sequence $\theta = (\theta_1, \theta_2, \ldots)$ of parameters by expanding $f$ on a fixed orthonormal basis $e_1, e_2, \ldots$ of $L_2(\mathcal{X}, \mathcal{A}, P_X)$. The Fourier coefficients take the form $\theta_i = \langle f, e_i \rangle = \mathrm{E} e_i(X) Y$.

The Fourier coefficients can be estimated unbiasedly by the estimators $Z_i = n^{-1} \sum_{j=1}^n Y_j e_i(X_j)$, but, as in Example 1.4, it appears not useful to try and reduce the model to the sequence model of Example 1.3 by considering $(Z_1, Z_2, \ldots)$ as the observation.

The assumption that the design distribution $P_X$ is known may be realistic in some practical situations, but is unpleasant. Perhaps it is a little surprising that it is not a merely technical assumption, but essential for the construction of our confidence sets. We intend to show elsewhere that the radius of the confidence sets will increase if $P_X$ is unknown, in varying amount, depending on what a priori assumptions are made on $P_X$. If $P_X$ is completely unknown, then intuitively this model should be equivalent to the fixed design regression model discussed in Example 1.6.

EXAMPLE 1.6 (Fixed regression).   In this model the observation is a vector $Y = (Y_1, \ldots, Y_n)$ of independent random variables distributed according to the regression model $Y_i = f(x_i) + \varepsilon_i$, for $\varepsilon_1, \ldots, \varepsilon_n$ i.i.d. normal variables with $\mathrm{E} \varepsilon_i = 0$ and $\mathrm{E} \varepsilon_i^2 = \sigma^2$ and $x_1, \ldots, x_n$ known constants. The variance $\sigma^2$ is known and the function $f$ is known to belong to a subset $\mathcal{F}$ of $L_2(\mathcal{X}, \mathcal{A}, \mu)$ for some distribution $\mu$.

Genovese and Wasserman [18] put this model in a sequence framework by expansion of the regression function on an empirical wavelet basis. They justify Beran [4] and Beran and Dümbgen [5] REACT confidence sets in terms of an honest confidence set over $\beta$-regular regression functions $f$, described in terms of a wavelet expansion. This is also the model treated by Juditsky and Lambert-Lacroix [25].

The model can be seen to reduce to a version of the finite sequence model of Example 1.2. All information about the regression function $f$ outside the design set $\{x_1, \ldots, x_n\}$ must stem from the model and not from the data. This point was made previously in Li [32], who gives the regression model as motivation for studying the finite sequence model. We shall not further discuss this model separately.

## 2. Construction of confidence regions.
Our method is based on sample splitting. We suppose that initial estimators $\hat{\theta}^{(n)}$ are given, and construct the confidence region based on $\hat{\theta}^{(n)}$ and an additional independent observation $X^{(n)}$. It was discussed previously how to split the data into independent "halves" that can be used for constructing $\hat{\theta}^{(n)}$ and $X^{(n)}$. The nature of the initial estimator $\hat{\theta}^{(n)}$ is irrelevant for the honesty of the confidence procedure, and hence $\hat{\theta}^{(n)}$ can be any of our favorite estimators. In particular, it can be an estimator that adapts to a selection of models of our choice.

Our confidence regions are based on estimators $R_n(\hat{\theta}^{(n)}) = R_n(\hat{\theta}^{(n)}, X^{(n)})$ of the squared norm $\|\theta - \hat{\theta}^{(n)}\|^2$ such that

$$(2.1) \quad \liminf_{n \to \infty} \inf_{\theta \in \Theta} P_\theta(R_n(\hat{\theta}^{(n)}) - \|\theta - \hat{\theta}^{(n)}\|^2 \geq -z_\alpha \hat{\tau}_{n,\theta} | \hat{\theta}^{(n)}) \geq 1 - \alpha,$$

for "scale estimators" $\hat{\tau}_{n,\theta}$ and "quantiles" $z_\alpha$. The probability is computed conditionally given the estimators $\hat{\theta}^{(n)}$, and hence refers only to the observation $X^{(n)}$ used to calculate $R_n(\hat{\theta}^{(n)})$ and $\hat{\tau}_{n,\theta}$. In view of Fatou's lemma the unconditional coverage probability will also be at least $1 - \alpha$. Then the set

$$(2.2) \qquad \hat{C}_n = \{\theta \in \Theta : \|\theta - \hat{\theta}^{(n)}\| \leq \sqrt{z_\alpha \hat{\tau}_{n,\theta} + R_n(\hat{\theta}^{(n)})}\}$$

is an honest confidence region with coverage probability at least $1 - \alpha$. (Define $\sqrt{x}$ to be 0 if $x < 0$.) The confidence region $\hat{C}_n$ is in general not a ball. However, in all our examples the scale estimators $\hat{\tau}_{n,\theta}$ satisfy

$$\hat{\tau}_{n,\theta} \lesssim \hat{\tau}_n + \frac{\|\theta - \hat{\theta}^{(n)}\|}{\sqrt{n}},$$

where $\lesssim$ denotes smaller than up to a constant which is fixed by the setting and $\hat{\tau}_n$ is independent of $\theta$ and determined by the size of the parameter set $\Theta$. It can be seen from this that the diameter of the confidence region satisfies

$$(2.3) \qquad \operatorname{diam}(\hat{C}_n) \lesssim \sqrt{\hat{\tau}_n} + \sqrt{R_n(\hat{\theta}^{(n)})} + n^{-1/2}.$$

(See the proof of the proposition below for a precise argument.) The last term on the right is the parametric rate of estimation and is typically negligible relative to the other terms. The first term $\sqrt{\hat{\tau}_n}$ depends on the supermodel $\Theta$ and its size is typically the same on every submodel.

The possibility of adaptation hinges on the second term. Typically (2.1) extends to a full, two-sided comparison, of the form $|R_n(\hat{\theta}^{(n)}) - \|\theta - \hat{\theta}^{(n)}\|^2| = O_{P_\theta}(\hat{\tau}_{n,\theta})$ uniformly in $\theta \in \Theta$. Then it follows that the diameter of $\hat{C}_n$ is of the order, uniformly in $\theta \in \Theta$,

$$\operatorname{diam}(\hat{C}_n) = O_{P_\theta}(\sqrt{\hat{\tau}_n} + \|\hat{\theta}^{(n)} - \theta\| + n^{-1/2}).$$

The diameter of the confidence set on a given submodel $\Theta_1 \subset \Theta$ is bounded above by the biggest order of the expression on the right-hand side under $\theta$, for $\theta$ ranging over $\Theta_1$. For small submodels, or more generally submodels where the estimators $\hat{\theta}^{(n)}$ perform well, the diameter will be dominated by the term $\sqrt{\hat{\tau}_n}$, the rate of the estimators of $\|\theta - \hat{\theta}^{(n)}\|^2$. On the other hand, in bigger submodels the term $\|\hat{\theta}^{(n)} - \theta\|$ may dominate. It is thus that we achieve adaptation to smaller models, but only up to the order $\sqrt{\hat{\tau}_n}$.

It is apparent from the preceding description that our confidence regions depend crucially on good estimators of the squared distance $\|\theta - \hat{\theta}^{(n)}\|^2$ of the parameter $\theta$ to the point $\hat{\theta}^{(n)}$. The latter point $\hat{\theta}^{(n)}$ may be considered fixed, as we condition on the initial estimator. The problem of constructing such estimators is therefore closely connected to the problem of estimating the squared norm $\|\theta\|^2$ of a Hilbert space-valued parameter. In some examples this is straightforward, but in the situations of density estimation and regression this problem is more involved. Fortunately, in the latter cases the estimation of a "quadratic functional" has been studied in detail by, among others, Fan [17], Bickel and Ritov [6], Laurent [26, 27] and Laurent and Massart [28], whose work obtains additional relevance in the present paper. The more recent papers consider adaptive estimators of the squared norm, but for our purposes optimal estimation under the biggest model will be sufficient. In view of their simplicity we shall adapt the constructions of Laurent [26, 27] to our purposes, but other approaches could be used as well.

This method consists of estimating the squared norm $\|\Pi_k\theta - \Pi_k\hat{\theta}^{(n)}\|^2$ of the projection of the difference $\theta - \hat{\theta}^{(n)}$ [where $\Pi_k\theta = (\theta_1, \ldots, \theta_k, 0, 0, \ldots)$] unbiasedly and trading off the resulting (squared) bias versus the variance of the estimator. Under the assumption that $\hat{\theta}^{(n)}$ takes its values in $\Theta$, the bias is bounded by a multiple of

$$(2.4) \qquad B_k^2 := \sup_{\theta \in \Theta} \|\theta - \Pi_k\theta\|^2.$$

The variance turns out to be of the order, for a parameter $\sigma^2$ that depends on the setting,

$$(2.5) \qquad \hat{\tau}_{k,n,\theta}^2 := \frac{2\sigma^4 k}{n^2} + \frac{4\sigma^2 \|\Pi_k\theta - \Pi_k\hat{\theta}^{(n)}\|^2}{n}.$$

The root $\hat{\tau}_{k,n,\theta}$ of this variance and the bias $B_k^2$ must be incorporated into the variable $\hat{\tau}_{n,\theta}$ as in (2.1). We define $\hat{\tau}_n = \sqrt{2}\sigma^2\sqrt{k}/n + B_k^2$, and conclude, in view of (2.3), that the diameter of the resulting confidence set (2.2) is of the order

$$(2.6) \qquad \frac{\sigma k^{1/4}}{\sqrt{n}} + B_k^2 + \|\theta - \hat{\theta}^{(n)}\| + \frac{\sigma}{\sqrt{n}}.$$

We can now choose an optimal value of $k$ by trading off $k^{1/4}/\sqrt{n}$ versus $B_k$.

The parameter $\sigma$ may depend on the unknown $\theta$, but in that case must be uniformly bounded over the supermodel $\Theta$.

For later reference we formalize the preceding as a proposition. Rather than making assumptions on bias and variance, we assume that the estimation rate of the estimators $R_{k,n}(\hat{\theta}^{(n)})$ is of the order as in the preceding discussion: for $\hat{\tau}_{k,n,\theta}$ as in (2.5), some number $z_\alpha$ and any sequences $k_n \to \infty$ and $M_n \to \infty$,

$$(2.7) \quad \limsup_{n \to \infty} \sup_{\theta \in \Theta} P_\theta(R_{k_n,n}(\hat{\theta}^{(n)}) - \|\Pi_{k_n}\theta - \Pi_{k_n}\hat{\theta}^{(n)}\|^2 \leq -z_\alpha \hat{\tau}_{k_n,n,\theta}) \leq \alpha,$$

$$(2.8) \quad \limsup_{n \to \infty} \sup_{\theta \in \Theta} P_\theta(|R_{k_n,n}(\hat{\theta}^{(n)}) - \|\Pi_{k_n}\theta - \Pi_{k_n}\hat{\theta}^{(n)}\|^2| \geq M_n \hat{\tau}_{k_n,n,\theta}) \to 0.$$

Of course, the second equation implies that the first is satisfied for sufficiently large $z_\alpha$, whereas an "absolute version" of the first equation for all $\alpha \in (0,1)$ will imply the second one.

PROPOSITION 2.1. *Suppose that $R_{k,n}(\hat{\theta}^{(n)})$ are estimators that satisfy* (2.7)–(2.8) *for $\hat{\tau}_{k,n,\theta}$ given in* (2.5) *and some $\sigma \in (0, \bar{\sigma}]$. Assume that $\hat{\theta}^{(n)}$ takes its values in $\Theta$. Then for $B_k$ given in* (2.4) *the sets*

$$\hat{C}_n = \{\theta \in \Theta : \|\theta - \hat{\theta}^{(n)}\| \leq \sqrt{z_\alpha \hat{\tau}_{k_n,n,\theta} + R_{k_n,n}(\hat{\theta}^{(n)})} + 2B_{k_n}\}$$

*are honest $(1-\alpha)$-confidence sets for $\theta \in \Theta$, for any $k_n \to \infty$, with diameter satisfying, for any $M_n \to \infty$,*

$$\limsup_{n \to \infty} \sup_{\theta \in \Theta} P_\theta\left(\operatorname{diam}(\hat{C}_n) \geq M_n\left[\frac{\bar{\sigma}k_n^{1/4}}{\sqrt{n}} + B_{k_n} + \|\theta - \hat{\theta}^{(n)}\|\right]\right) = 0.$$

PROOF. By (2.4) the difference $\|\theta - \hat{\theta}^{(n)}\|$ is bounded above by $\|\Pi_k(\theta - \hat{\theta}^{(n)})\| + 2B_k$. Therefore, by the definition of $\hat{C}_n$,

$$P_\theta(\theta \notin \hat{C}_n) \leq P_\theta(\|\Pi_k(\theta - \hat{\theta}^{(n)})\| > \sqrt{z_\alpha \hat{\tau}_{k_n,n,\theta} + R_{k_n,n}(\hat{\theta}^{(n)})}).$$

In view of (2.7) the right-hand side is asymptotically bounded above by $\alpha$, uniformly in $\theta \in \Theta$. Hence $\hat{C}_n$ is an asymptotic confidence region of confidence level $1 - \alpha$.

In view of the form (2.5) of $\hat{\tau}_{k,n,\theta}$ every element $\theta$ of $\hat{C}_n$ satisfies

$$\|\theta - \hat{\theta}^{(n)}\| \leq \sqrt{z_\alpha \frac{\sqrt{2}\sigma^2\sqrt{k}}{n} + R_{k_n,n}(\hat{\theta}^{(n)})} + 2B_{k_n} + \frac{\sqrt{z_\alpha 2\sigma}}{n^{1/4}}\sqrt{\|\theta - \hat{\theta}^{(n)}\|}.$$

The inequality $x \leq B + A\sqrt{x}$ for real numbers $x$ and positive real numbers $A$ and $B$ implies that $x \leq 2B + 2A^2$. We conclude that the diameter of $\hat{C}_n$ is bounded by a multiple of

$$\frac{\bar{\sigma}k^{1/4}}{\sqrt{n}} + \sqrt{R_{k_n,n}(\hat{\theta}^{(n)})} + B_{k_n} + \frac{\bar{\sigma}}{\sqrt{n}}.$$

The variables $R_{k_n,n}(\hat{\theta}^{(n)})$ are with $P_\theta$-probability tending to 1 bounded above by a multiple of $\|\Pi_{k_n}\theta - \Pi_{k_n}\hat{\theta}^{(n)}\|^2 + M_n\hat{\tau}_{k_n,n,\theta}$, for any given $M_n \to \infty$, by (2.8). Therefore, with probability tending to 1 the diameter of $\hat{C}_n$ is bounded by

$$\frac{\bar{\sigma} k^{1/4}}{\sqrt{n}} + \|\theta - \hat{\theta}^{(n)}\| + \sqrt{M_n \hat{\tau}_{k,n,\theta}} + B_{k_n} + \frac{\bar{\sigma}}{\sqrt{n}}.$$

Here the last term is negligible relative to the first. The proposition follows in view of the form (2.5) of $\hat{\tau}_{k,n,\theta}$ and the inequality $\sigma k^{1/4}/\sqrt{n} + \sqrt{\sigma}\sqrt{x}/n^{1/4} + x \le 2\sigma k^{1/4}/\sqrt{n} + 2x$, which is valid for any $k \ge 1$, $x \ge 0$ and $\sigma > 0$. □

The natural (or "naive") estimators $R_n(\hat{\theta}^{(n)})$ of $\|\theta - \hat{\theta}^{(n)}\|^2$ in our examples assume negative values, which could lead to a confidence set $\hat{C}_n$ in (2.2) of zero diameter. This is unlike the usual situation in parametric models, where $\sqrt{n}$ times the radius of a confidence region for $\theta$ generally has the desirable property of tending in probability to a positive constant. In practice it might be useful to eliminate the possibility of radii of zero by substituting for the right-hand side of (2.2) a more conservative cut-off, given by the maximum of the current right-hand side and $\sqrt{z_\alpha \hat{\tau}_{n,\theta}}$ (or perhaps $\sqrt{z_\alpha \hat{\tau}_{n,\theta}/2}$).

EXAMPLE 2.1 (Model of dimension $n$).   If $\Theta = \mathbb{R}^n$, then we can avoid a bias by choosing $k = n$. Then the diameter of the confidence sets is of the order equal to the maximum of $n^{-1/4}$ and the estimation error $\|\theta - \hat{\theta}^{(n)}\|$.

EXAMPLE 2.2 (Regular models).   The usual models to define regular parameter are the ellipsoids $S(\beta, L) = \{\theta \in \ell_2 : \sum_{i=1}^{\infty} \theta_i^2 i^{2\beta} \le L^2\}$, for $\beta > 0$ and $L > 0$ given. Suppose we choose $\Theta = S(\beta, L)$ for fixed values of $\beta$ and $L$ as the supermodel, on which we require honesty, and consider adaptation on ellipsoids defined by different parameter values.

If we cut off the series expansion at level $k$, then the maximal squared bias is equal to

$$\sup_{\theta \in S(\beta,L)} \sum_{i=k+1}^{\infty} \theta_i^2 \le \sup_{\theta \in S(\beta,L)} \sum_{i=k+1}^{\infty} \theta_i^2 \left(\frac{i}{k}\right)^{2\beta} \le \frac{L^2}{k^{2\beta}}.$$

This leads to the trade-off $k^{1/4}/\sqrt{n} \sim L/k^\beta$, resulting in a cut-off of the order

$$k \sim L^{4/(4\beta+1)} n^{1/(2\beta+1/2)}$$

and a bias of the order $n^{-\beta/(2\beta+1/2)} L^{1/(4\beta+1)}$.

This choice of $k$ is compatible with $k \le n$ only if $\beta \ge 1/4$. Thus if $\theta$ is restricted to $\mathbb{R}^n$, as in the finite sequence model, then we consider submodels $S(\beta, L)$ with $\beta \ge 1/4$ only.

For this choice of $k$ we obtain a confidence region for the full parameter $\theta \in \ell_2$ of diameter of order equal to the maximum of order $n^{-\beta/(2\beta+1/2)}$ and the estimation error $\|\theta - \hat{\theta}^{(n)}\|$. The lower bound $n^{-\beta/(2\beta+1/2)}$ and the cut-off $k$ are for some values of $\beta$ given in the third and fourth columns of Table 1.

Thus the role of the minimal diameter $n^{-1/4}$ in the preceding example is now taken over by $n^{-\beta/(2\beta+1/2)}$.

For the initial estimators $\hat{\theta}^{(n)}$ there is a variety of choices. A relatively simple scheme is to choose $\hat{\theta}^{(n)}$ to adapt to all regularity classes $S(\gamma, M)$ in the sense that, for all $\gamma \geq \beta$ and all $M > 0$, for some constants $C_{\gamma,M}$,

$$\sup_{\theta \in S(\gamma, M)} \mathrm{E}_\theta (\hat{\theta}^{(n)} - \theta)^2 \leq C_{\gamma,M} n^{-2\gamma/(2\gamma+1)}.$$

Such estimators exist in the examples considered in the Introduction. In fact, there exist estimators that adapt to a much larger collection of submodels than only the Sobolev models considered in this paper. Combined with our construction this will lead to a confidence region around $\hat{\theta}^{(n)}$ of diameter of the order $n^{-\gamma/(2\gamma+1)}$ uniformly over $S(\gamma, M)$ if $\gamma \in [\beta, 2\beta]$, and of the order $n^{-\beta/(2\beta+1/2)}$ over $S(\gamma, M)$ for other indices $\gamma$.

**3. Sequence models.** Suppose that we observe an infinite sequence $X = (X_1, X_2, \ldots)$ of independent random variables $X_i$ possessing means $\mathrm{E}X_i = \theta_i$ and variances $\sigma^2/n$. The parameter $\theta = (\theta_1, \theta_2, \ldots)$ is known to belong to a subset $\Theta$ of $\ell_2$. This formulation encompasses both the finite and the infinite sequence models of Examples 1.2 and 1.3, if in the former case it is understood that $\Theta \subset \mathcal{R}_n := \{\theta \in \ell_2 : \theta_i = 0, i > n\}$ and that $X_{n+1}, X_{n+2}, \ldots$ may not be used to estimate $\sigma^2$. Our main interest is in the case where the $X_i$ are also normally distributed, but we also consider the more general situation. The assumption of normality allows a precise and simple derivation of the radius of a confidence region. In a final subsection we also indicate how to obtain confidence sets with guaranteed level for finite $n$.

3.1. *Normal distributions.* In this section we assume in addition to the preceding that each $X_i$ is normally distributed.

Given an initial estimator $\hat{\theta}^{(n)}$, based on observations that are independent of $X$, our estimator for $\|\theta - \hat{\theta}^{(n)}\|^2$ is given by

$$(3.1) \qquad R_{k,n}(\hat{\theta}^{(n)}) = \sum_{i=1}^{k} (X_i - \hat{\theta}_i^{(n)})^2 - \frac{k\sigma^2}{n}.$$

Here $k = k_n$ is chosen dependent on $\Theta$ and/or $\Theta_1$, where we must have $k \leq n$ in the finite sequence model. This estimator is combined with the estimator of variance (random only in its dependence on $\hat{\theta}^{(n)}$)

$$(3.2) \qquad \hat{\tau}_{k,n,\theta}^2 = \frac{2k\sigma^4}{n^2} + \frac{4\sigma^2}{n} \sum_{i=1}^{k} (\theta_i - \hat{\theta}_i^{(n)})^2.$$

We shall show that $R_{k,n}(\hat{\theta}^{(n)})$ tends in distribution to a normal distribution, uniformly in $\theta \in \ell_2$. This allows us to construct confidence sets of the type (2.2) by using normal quantiles for the values $z_\alpha$. [Because $R_{k,n}$ and $\hat{\tau}_{k,n,\theta}$ depend in fact only on $(\theta_1, \ldots, \theta_k)$, "uniformly in $\theta \in \ell_2$" means effectively "uniformly in $(\theta_1, \ldots, \theta_k) \in \mathbb{R}^k$."] Because $R_{k,n}(\hat{\theta}^{(n)})$ is a sum of independent variables, its asymptotic normality is not a surprise. The main contribution of the following theorem is that this asymptotic normality is uniform in $\theta$, without any conditions on the initial estimators $\hat{\theta}^{(n)}$. This depends crucially on the normality of the observations.

The convergence in the following theorem may be understood in the almost sure sense. As the proof shows, the weak convergence is actually uniform in the values $\hat{\theta}^{(n)}$.

THEOREM 3.1. *For any $k_n \to \infty$ as $n \to \infty$,*

$$\sup_{\theta \in \ell_2} \sup_{x \in \mathbb{R}} \left| P_\theta \left( \frac{R_{k_n,n}(\hat{\theta}^{(n)}) - \sum_{i=1}^{k_n}(\theta_i - \hat{\theta}_i^{(n)})^2}{\hat{\tau}_{k_n,n,\theta}} \le x \Big| \hat{\theta}^{(n)} \right) - \Phi(x) \right| \to 0.$$

PROOF. We can express the variable $(R_{k,n}(\hat{\theta}^{(n)}) - \sum_{i=1}^{k}(\theta_i - \hat{\theta}_i^{(n)})^2)/\hat{\tau}_{k,n,\theta}$ in the independent standard normal variables $\varepsilon_i$ defined by $X_i = \theta_i + (\sigma/\sqrt{n})\varepsilon_i$ as

$$\frac{1}{\hat{\tau}_{k,n,\theta}} \left( \sum_{i=1}^{k} (\varepsilon_i^2 - 1)\frac{\sigma^2}{n} + \frac{2\sigma}{\sqrt{n}} \sum_{i=1}^{k}(\theta_i - \hat{\theta}_i^{(n)})\varepsilon_i \right)$$

$$= \frac{1}{\sqrt{2k}} \sum_{i=1}^{k}(\varepsilon_i^2 - 1)A_{n,k}(\theta) + \frac{\sum_{i=1}^{k}(\theta_i - \hat{\theta}_i^{(n)})\varepsilon_i}{\sqrt{\sum_{i=1}^{k}(\theta_i - \hat{\theta}_i^{(n)})^2}} B_{n,k}(\theta),$$

for the positive constants whose squares are given by

$$A_{k,n}(\theta)^2 = \frac{1}{1 + (2n/k\sigma^2)\sum_{i=1}^{k}(\theta_i - \hat{\theta}_i^{(n)})^2},$$

$$B_{k,n}(\theta)^2 = \frac{\sum_{i=1}^{k}(\theta_i - \hat{\theta}_i^{(n)})^2}{(k\sigma^2/2n) + \sum_{i=1}^{k}(\theta_i - \hat{\theta}_i^{(n)})^2}.$$

By the rotational invariance of the multivariate standard normal distribution, for any vector $\psi$ with norm 1 the random vector $((2k)^{-1/2}\sum_{i=1}^{k}(\varepsilon_i^2 - 1), \sum_{i=1}^{k}\psi_i\varepsilon_i)$ is equal in distribution to the random vector $((2k)^{-1/2}\sum_{i=1}^{k}(\varepsilon_i^2 - 1), k^{-1/2}\sum_{i=1}^{k}\varepsilon_i)$. The latter vector tends in distribution to a vector of two independent standard normal variables, as $n \to \infty$. The coefficients $A_{k,n}(\theta)$ and $B_{k,n}(\theta)$ are contained in the unit interval and satisfy $A_{k,n}^2(\theta) + B_{k,n}^2(\theta) = 1$, for any $k, n, \theta$.

We can complete the proof by noting that if a sequence of random vectors $(X_n, Y_n)$ converges in distribution to a random vector $(X, Y)$, then the sequence $AX_n + BY_n$ tends in distribution to $AX + BY$, uniformly in coefficients $(A, B)$ belonging to a compact set. $\square$

The theorem shows that $R_{k,n}(\hat{\theta}^{(n)})$ is a good estimator of the squared norm of the projection $\Pi_k(\theta - \hat{\theta}^{(n)})$ of $\theta - \hat{\theta}^{(n)}$ onto the $k$-dimensional subspace $\{\theta \in \ell_2 : \theta_i = 0, i > k\}$, and justifies (2.8) with $\hat{\tau}^2_{k,n,\theta}$ of the order as in (2.5). Thus Proposition 2.1 yields a confidence region of diameter of the order

$$\bar{\sigma}\left(\frac{k_n}{n^2}\right)^{1/4} + B_{k_n} + \|\theta - \hat{\theta}^{(n)}\|.$$

EXAMPLE 3.1 (Finite sequence model). In the finite sequence model of Example 1.2 with $\Theta = \mathbb{R}^n$, we have bias $B_k$ zero if we choose $k = n$. This leads to confidence sets of diameter of the order equal to the maximum of $n^{-1/4}$ and the estimation error $\|\theta - \hat{\theta}^{(n)}\|$.

As was shown by Li [32] and Baraud [1] the $n^{-1/4}$ lower bound cannot be improved upon without losing full honesty.

We can influence the term $\|\theta - \hat{\theta}^{(n)}\|$ by choosing our favorite estimators $\hat{\theta}^{(n)}$. For instance, we may choose any of the adaptive penalized minimum contrast estimators considered in [8]. As shown by Birgé and Massart [8] we can adapt to large classes of a priori models by choosing appropriate penalties.

One choice of penalties leads to estimators that, among other good properties, satisfy, for every $D$,

$$\sup_{\theta \in \Theta_D} \mathrm{E}_\theta \|\hat{\theta}^{(n)} - \theta\|^2 \lesssim \frac{\sigma^2}{n}\left[D + \log\left(\frac{2n}{D}\right) + 1\right],$$

where $\Theta_D = \{\theta \in \mathbb{R}_n : \#(\theta_i \neq 0) \leq D\}$. The confidence sets centered at these estimators attain a uniform order equal to the maximum of $n^{-1/4}$ and $\sqrt{D/n} + \sqrt{\log(2n/D)/n}$. As long as $D \ll n$ this improves upon the order 1 rate attained by the naive chi-square procedure, and we obtain the best possible rate $n^{-1/4}$ uniformly over every set $\Theta_D$ with $D \leq \sqrt{n}$. Thus these excellent adaptation properties of $\hat{\theta}^{(n)}$ result in smaller confidence regions, for more submodels, than those found in [1], pages 533–536, by a direct construction.

The estimators $R_{k,n}(\hat{\theta}^{(n)})$ and $\hat{\tau}_{k,n,\theta}$ in the preceding theorem depend on $\sigma^2$ and hence so far we have implicitly assumed that (an upper bound on) the variance $\sigma^2$ is known. The preceding remains true if it is replaced by a good estimator.

THEOREM 3.2. *The assertion of Theorem* 3.1 *remains true if* $\sigma^2$ *in the definitions* (3.1) *and* (3.2) *of* $R_{k,n}$ *and* $\hat{\tau}_{k,n,\theta}$ *is replaced by estimators* $\hat{\sigma}_n^2$ *such that*

$$\sup_{\theta \in \Theta} \mathrm{P}_\theta (\sqrt{k}_n |\hat{\sigma}_n^2 - \sigma^2| > \varepsilon \mid \hat{\theta}^{(n)}) \to 0.$$

PROOF.    Represent the observations as $X_i = \theta_i + (\sigma/\sqrt{n})e_i$ for independent standard normal variables $e_i$. It suffices to prove the uniform asymptotic normality of the variables

$$\frac{\sum_{i=1}^k (e_i^2 - 1)\sigma^2/n + k(\sigma^2/n - \hat{\sigma}^2/n) + (2\sigma/\sqrt{n})\sum_{i=1}^k (\theta_i - \hat{\theta}_i^{(n)})e_i}{\sqrt{(2\hat{\sigma}^4 k/n^2) + (4\hat{\sigma}^2/n)\sum_{i=1}^k (\theta_i - \hat{\theta}_i^{(n)})^2}}.$$

Therefore, it suffices to prove that, uniformly in $\theta \in \Theta$,

(3.3)    $$\frac{2\hat{\sigma}^4/n^2 + 4\hat{\sigma}^2/nk \sum_{i=1}^k (\theta_i - \hat{\theta}_i^{(n)})^2}{2\sigma^4/n^2 + (4\sigma^2/nk)\sum_{i=1}^k (\theta_i - \hat{\theta}_i^{(n)})^2} - 1 \xrightarrow{\mathrm{P}} 0,$$

(3.4)    $$\frac{\sqrt{k}}{n} \frac{(\hat{\sigma}^2 - \sigma^2)}{\sqrt{2\sigma^4/n^2 + (4\sigma^2/nk)\sum_{i=1}^k (\theta_i - \hat{\theta}_i^{(n)})^2}} \xrightarrow{\mathrm{P}} 0.$$

The absolute value of the left-hand side of (3.3) can be rewritten in the form, for the constants $C_{n,k}(\theta) = 2n/(k\sigma^2)\sum_{i=1}^k (\theta_i - \hat{\theta}_i^{(n)})^2$,

$$\frac{\hat{\sigma}^2}{\sigma^2} \left| \frac{\hat{\sigma}^2/\sigma^2 + C_{n,k}(\theta)}{1 + C_{n,k}(\theta)} - 1 \right| \leq \left| \frac{\hat{\sigma}^2}{\sigma^2} - 1 \right|.$$

Thus this reduces to $\hat{\sigma}/\sigma \xrightarrow{\mathrm{P}} 1$, uniformly in $\theta$. Assertion (3.4) is true as soon as $\sqrt{k}(\hat{\sigma}^2 - \sigma^2) \xrightarrow{\mathrm{P}} 0$, uniformly in $\theta$.    □

In the finite sequence model with $\Theta = \mathbb{R}^n$ there is no possibility to estimate $\sigma^2$, and the same is true in the infinite sequence model without some restriction on the parameter set $\Theta$. On the other hand, in the infinite sequence model with a restriction to regular parameters, estimation of $\sigma^2$ is easy.

EXAMPLE 3.2 (Regular models).    For given integers $l$ and $m$ consider the estimator $\hat{\sigma}^2 = (n/l)\sum_{i=m+1}^{m+l} X_i^2$. This has mean and variance given by

$$\mathrm{E}_\theta \hat{\sigma}^2 = \frac{n}{l} \sum_{i=m+1}^{m+l} \left( \frac{\sigma^2}{n} + \theta_i^2 \right) = \sigma^2 + \frac{n}{l}\sum_{i=m+1}^{m+l} \theta_i^2,$$

$$\mathrm{var}_\theta \hat{\sigma}^2 = \frac{n^2}{l^2} \sum_{i=m+1}^{m+l} \left( \frac{4\sigma^2\theta_i^2}{n} + \frac{2\sigma^4}{n^2} \right) = \left( \frac{4n\sigma^2}{l^2} \sum_{i=m+1}^{m+l} \theta_i^2 + \frac{2\sigma^4}{l} \right).$$

It follows that the mean squared error over the regularity class $S(\beta, L)$ can be bounded as

$$\sup_{\theta \in S(\beta, L)} \mathrm{E}_{\theta}(\hat{\sigma}^2 - \sigma^2)^2 \lesssim \frac{n^2}{l^2} \frac{1}{m^{4\beta}} + \frac{n}{l^2} \frac{1}{m^{2\beta}} + \frac{1}{l}.$$

In view of Theorem 3.2 we wish this to be of smaller order than $1/k$.

In the infinite sequence model with $\Theta = S(\beta, L)$ as the biggest model, we choose $k = n^{1/(2\beta+1/2)}$ (cf. Example 2.2), and hence we must choose $l \gg n^{1/(2\beta+1/2)}$. These values are shown for some values of $\beta$ in Table 1. For the minimal value of $l$ we must choose $m \geq n^{1/(2\beta+1/2)}$ and then the estimator for $\sigma^2$ becomes independent of $R_{k,n}(\hat{\theta}^{(n)})$. A variety of other combinations of $(m, l)$ will do as well.

In the finite sequence model with $\Theta$ restricted to $S(\beta, L)$, truncated to $\mathbb{R}^n$, the choice $l \gg n^{1/(2\beta+1/2)}$ can be realized with $l \leq n$ only if $\beta > 1/4$. We can then combine it with $m$ of the order $n^{1/(2\beta+1/2)}$.

### 3.2. *Nonnormal distributions.*

The assumed normality of the observations $X_1, X_2, \ldots$ in the preceding section helps one to obtain precise critical values, but it is not important for the general ideas. In this section let $X_i = \theta_i + (\sigma/\sqrt{n})\varepsilon_i$ for an i.i.d. sequence $\varepsilon_1, \varepsilon_2, \ldots$ with mean zero, variance 1 and finite fourth moment. Then define $R_{k,n}(\hat{\theta}^{(n)})$ as in (3.1) and define the variance estimator

$$\hat{\tau}_{k,n,\theta}^2 = \frac{k\sigma^4 \operatorname{var}(\varepsilon_1^2)}{n^2} + \frac{4\sigma^2}{n} \sum_{i=1}^{k} (\theta_i - \hat{\theta}_i^{(n)})^2 + \frac{4\sigma^3}{n\sqrt{n}} \sum_{i=1}^{k} (\theta_i - \hat{\theta}_i^{(n)}) \operatorname{cov}(\varepsilon_1^2, \varepsilon_1).$$

THEOREM 3.3. *For any $k$ and $n$,*

$$\inf_{\theta \in \ell_2} \mathrm{P}_{\theta}\left( \left| \frac{R_{k,n}(\hat{\theta}^{(n)}) - \sum_{i=1}^{k}(\theta_i - \hat{\theta}_i^{(n)})^2}{\hat{\tau}_{k,n,\theta}} \right| \leq \frac{1}{\sqrt{\alpha}} \Big| \hat{\theta}^{(n)} \right) \geq 1 - \alpha.$$

PROOF. The quantity $R_{k,n}(\hat{\theta}^{(n)})$ is an unbiased estimator of $\sum_{i=1}^{k}(\theta_i - \hat{\theta}_i^{(n)})^2$, and $\hat{\tau}_{k,n,\theta}^2$ is equal to its variance. Therefore, the inequality follows by Chebyshev's inequality. □

The preceding theorem is based on Chebyshev's inequality, which is notably imprecise. However, this crude device costs only in terms of the constants and not in terms of the rate. If $Z$ is exactly standard normal distributed, then we have that $\mathrm{P}(|Z| \geq 1.96) = 0.05$, whereas the use of Chebyshev's inequality $\mathrm{P}(|Z| \geq M) \leq M^{-2}$ would replace the normal quantile 1.96 by $M = 0.05^{-1/2} \approx 4.5$, so that the resulting confidence set would be a bit more than two times too wide.

For many estimators $\hat{\theta}^{(n)}$ we can avoid this penalty, because the quantities $R_{k,n}(\hat{\theta}^{(n)})$ will be asymptotically normal, at least under the overall probability law governing the initial estimators $\hat{\theta}^{(n)}$ and the observations $X^{(n)}$. This will depend on the initial estimators $\hat{\theta}^{(n)}$, but the following assumption appears to be reasonable. Assume that the initial estimators satisfy, for some sequence $\varepsilon_n \to 0$,

$$(3.5) \qquad \sup_{\theta \in \Theta} P_\theta \left( \max_{1 \le i \le k_n} |\hat{\theta}_i^{(n)} - \theta_i|^2 > \varepsilon_n \sum_{i=1}^{k_n} (\hat{\theta}_i^{(n)} - \theta_i)^2 \right) \to 0.$$

THEOREM 3.4.  *For any $k_n \to \infty$ as $n \to \infty$ such that* (3.5) *holds,*

$$\sup_{\theta \in \ell_2} \sup_{x \in \mathbb{R}} \left| P_\theta \left( \frac{R_{k_n,n}(\hat{\theta}^{(n)}) - \sum_{i=1}^{k_n} (\theta_i - \hat{\theta}_i^{(n)})^2}{\hat{\tau}_{k_n,n,\theta}} \le x \right) - \Phi(x) \right| \to 0.$$

PROOF.  We can express the variable $(R_{k,n}(\hat{\theta}^{(n)}) - \sum_{i=1}^{k} (\theta_i - \hat{\theta}_i^{(n)})^2)/\hat{\tau}_{k,n,\theta}$ in the form

$$(3.6) \qquad \sum_{i=1}^{k} A_{k,n}(\theta)(\varepsilon_i^2 - 1) + \sum_{i=1}^{k} B_{i,k,n}(\theta)\varepsilon_i,$$

for the positive constants given by

$$A_{k,n}(\theta) = \frac{\sigma^2}{n\hat{\tau}_{k,n,\theta}},$$

$$B_{k,n}(\theta) = \frac{2\sigma(\theta_i - \hat{\theta}_i^{(n)})}{\sqrt{n}\hat{\tau}_{k,n,\theta}}.$$

The terms in the sum (3.6) are conditionally independent under $P_\theta^{(n)}$ given $\hat{\theta}^{(n)}$, and the sum has conditional mean and variance equal to 0 and 1. If the terms of the sum also satisfy the conditional Lindeberg condition in probability, then the variables (3.6) converge conditionally in distribution, in probability. We wish to show that this is true uniformly in $\theta \in \Theta$.

Thus it suffices to prove that for every $k_n \to \infty$, every $\delta > 0$ and any sequence $\{\theta_n\} \subset \Theta$ as $n \to \infty$,

$$\sum_{i=1}^{k_n} E_{\theta_n}(A_{k_n,n}(\theta_n)(\varepsilon_i^2 - 1) + B_{i,k_n,n}(\theta_n)\varepsilon_i)^2 \mathbb{1}_{|A_{k_n,n}(\theta_n)(\varepsilon_i^2-1)+B_{i,k_n,n}(\theta_n)\varepsilon_i| > \delta} \to 0.$$

For any $c \in [0,1)$ and positive numbers $A, B$ we have that $(1-c)(A^2 + B^2) \le A^2 + B^2 - 2cAB$. Because the correlation $c$ between $\varepsilon_1^2$ and $\varepsilon_1$ is nonnegative and strictly smaller than 1, this inequality can be used to see that

$$(1-c)\left( \frac{k\sigma^4 \operatorname{var}(\varepsilon_1^2)}{n^2} + \frac{4\sigma^2}{n} \sum_{i=1}^{k} (\theta_i - \hat{\theta}^{(n)})^2 \right) \le \hat{\tau}_{k,n,\theta}^2.$$

Consequently,

$$\max_{1 \leq i \leq k}(A^2_{k,n}(\theta) + B^2_{i,k,n}(\theta)) \lesssim \frac{1}{k} + \frac{\max_{1 \leq i \leq k}(\theta_i - \hat{\theta}^{(n)}_i)^2}{\sum_{i=1}^k(\theta_i - \hat{\theta}^{(n)}_i)^2}.$$

By assumption the right-hand side converges to zero in probability, as $k = k_n \to \infty$ and $n \to \infty$. We also have that $\sum_{i=1}^k(A^2_{k,n}(\theta) + B^2_{i,k,n}(\theta))$ is uniformly bounded. We can conclude that the Lindeberg condition is satisfied. □

3.3. *Exact simulation.* The procedures in the preceding section can be implemented as soon as the lower-order moments of the errors $\varepsilon_i$ are known (or can be estimated). If the full distribution of the errors is available, then we may also obtain exact, finite-sample confidence regions. This observation is even of interest in the case of Gaussian errors.

The variable $(R_{k,n}(\hat{\theta}^{(n)}) - \sum_{i=1}^k(\theta_i - \hat{\theta}^{(n)}_i)^2)/\hat{\tau}_{k,n,\theta}$ can be written as a function $S_{k,n}(\varepsilon_1, \ldots, \varepsilon_n, \theta, \hat{\theta}^{(n)})$, as in (3.6) in the proof of Theorem 3.4. This representation allows simulation of the distribution of the given variable under $\theta$, for every fixed $\theta \in \Theta$. Thus in principle we can find the $\alpha$-quantile $-z_\alpha(\theta)$ of this distribution, for every $\theta$. Then $\hat{C}_n$ given in Proposition 2.1, but with $z_\alpha$ replaced by $z_\alpha(\theta)$, is a valid $(1 - \alpha)$-confidence region.

Under the conditions of Theorem 3.4 the quantiles $z_\alpha(\theta)$ converge to Gaussian quantiles, uniformly in $\theta$.

## 4. Density estimation.
Suppose that we observe an i.i.d. sample $X_1, \ldots, X_n$ from a density $f$ relative to some measure $\mu$ on a measurable space $(\mathcal{X}, \mathcal{A})$. Let $\theta = (\theta_1, \theta_2, \ldots)$ be the Fourier coefficients of $f$ relative to a given orthonormal basis of $L_2(\mathcal{X}, \mathcal{A}, \mu)$, and let $\Theta$ correspond to the collection of all densities deemed possible. Assume that the densities $\theta \in \Theta$ are uniformly bounded.

Given an initial estimator $\hat{\theta}^{(n)}$ our estimator for $\|\theta - \hat{\theta}^{(n)}\|^2$ is given by

$$R_{k,n}(\hat{\theta}^{(n)}) = \frac{1}{n(n-1)} \sum_{r \neq s=1}^n \sum_{i=1}^k (e_i(X_r) - \hat{\theta}^{(n)}_i)(e_i(X_s) - \hat{\theta}^{(n)}_i).$$

Here $k = k_n$ is chosen dependent on $\Theta$. We combine this with the variance estimator

$$\hat{\tau}^2_{k,n,\theta} = \frac{2k\|f\|^2_\infty}{n(n-1)} + \frac{4\|f\|_\infty}{n} \sum_{i=1}^k(\theta_i - \hat{\theta}^{(n)}_i)^2.$$

THEOREM 4.1. *For any $k, n$,*

$$\sup_{\theta \in \Theta} \mathrm{E}_\theta\left(\left(\frac{R_{k,n}(\hat{\theta}^{(n)}) - \sum_{i=1}^k(\theta_i - \hat{\theta}^{(n)}_i)^2}{\hat{\tau}_{k,n,\theta}}\right)^2 \Big| \hat{\theta}^{(n)}\right) \leq 1.$$

PROOF.   The estimator $R_{k,n}(\hat{\theta}^{(n)})$ is a $U$-statistic of order 2 with kernel $h(x,y) = \sum_{i=1}^k (e_i(x) - \hat{\theta}_i^{(n)})(e_i(y) - \hat{\theta}_i^{(n)})$. Its mean is equal to

$$\mathrm{E}h(X_1, X_2) = \sum_{i=1}^k (\theta_i - \hat{\theta}_i^{(n)})^2.$$

Its Hoeffding decomposition (e.g., [39], Section 11.4) is

$$
\begin{aligned}
R_{k,n}(\hat{\theta}^{(n)}) = {} & \mathrm{E}h(X_1, X_2) + \frac{1}{n}\sum_{r=1}^n P_1 h(X_r) \\
& + \frac{1}{n(n-1)}\sum_{r \neq s}\sum P_{1,2}h(X_r, X_s),
\end{aligned}
$$

(4.1)

for the "kernel functions" given by

$$P_1 h(x) = 2\sum_{i=1}^k (\theta_i - \hat{\theta}_i^{(n)})(e_i(x) - \theta_i),$$

$$
\begin{aligned}
P_{1,2}h(x,y) &= \sum_{i=1}^k (e_i(x) - \theta_i)(e_i(y) - \theta_i) \\
&= \sum_{i=1}^k e_i(x)e_i(y) - \sum_{i=1}^k \theta_i(e_i(x) + e_i(y)) + \sum_{i=1}^k \theta_i^2.
\end{aligned}
$$

The three terms of the Hoeffding decomposition and also each of the individual terms in its sums are uncorrelated. Furthermore, the variance of the last term in (4.1) is equal to $2/(n(n-1))\,\mathrm{var}\,P_{1,2}h(X_1, X_2)$.

The variance of a factor in the linear term can be bounded as

$$
\begin{aligned}
\mathrm{var}(P_1 h(X_1)) &= 4\mathrm{E}\left(\sum_{i=1}^k (\theta_i - \hat{\theta}_i^{(n)})e_i(X_1)\right)^2 \\
&\leq 4\|f\|_\infty \int \left(\sum_{i=1}^k (\theta_i - \hat{\theta}_i^{(n)})e_i\right)^2 d\mu = 4\|f\|_\infty \sum_{i=1}^k (\theta_i - \hat{\theta}_i^{(n)})^2,
\end{aligned}
$$

by the orthonormality of the functions $e_i$ in $L_2(\mu)$.

The variables $\sum_{i=1}^k (e_i(X_1) - \theta_i)(e_i(X_2) - \theta_i)$ and $\sum_{i=1}^k \theta_i(e_i(X_1) + e_i(X_2))$ are uncorrelated and their sum is $\sum_{i=1}^k e_i(X_1)e_i(X_2) + \sum_{i=1}^k \theta_i^2$. It follows that

$$\mathrm{var}(P_{1,2}h(X_1, X_2)) = \mathrm{var}\sum_{i=1}^k e_i(X_1)e_i(X_2) - \mathrm{var}\sum_{i=1}^k \theta_i(e_i(X_1) + e_i(X_2)).$$

This becomes bigger if we leave out the second variance on the right and replace the first variance on the right by the second moment $E(\sum_{i=1}^{k} e_i(X_1) \times e_i(X_2))^2$, which can be bounded by

$$\|f\|_\infty^2 \iint \left(\sum_{i=1}^{k} e_i(x) e_i(y)\right)^2 d\mu(x)\,d\mu(y) = k\|f\|_\infty^2,$$

by the orthonormality of the functions $e_i$ in $L_2(\mu)$. $\square$

By Markov's inequality, if we choose $z_\alpha = \sqrt{1/\alpha}$, then

$$\inf_{\theta \in \Theta} P_\theta\left(\left|R_{k,n}(\hat{\theta}^{(n)}) - \sum_{i=1}^{k}(\theta_i - \hat{\theta}_i^{(n)})^2\right| \leq z_\alpha \hat{\tau}_{k,n,\theta} \mid \hat{\theta}^{(n)}\right) \geq 1 - \alpha.$$

The present variance $\hat{\tau}_{k,n,\theta}^2$ has exactly the same form as in Section 3, with $\|f\|_\infty$ playing the role of $\sigma^2$. For more precision we can express $\|f\|_\infty$ in $\theta$, and it is not necessary to know a uniform bound on the regression functions. The approximation (2.8) with $\hat{\tau}_{k,n,\theta}^2$ of the order as in (2.5) is again satisfied and Proposition 2.1 yields a confidence region of diameter of the order, with $M$ a uniform bound on $\Theta$,

$$M\left(\frac{k_n}{n^2}\right)^{1/4} + B_{k_n} + \|\theta - \hat{\theta}^{(n)}\|.$$

The corollaries for, for example, regular models are the same.

Depending on the basis functions $e_i$, the resulting confidence region can be tightened by using higher moments or exponential bounds. Finding an exact limit distribution appears to be not straightforward. Existing limit results for $U$-statistics with changing kernels (e.g. [33]) are based on approximation of the kernel by a finite product kernel of fixed dimension. In our case the kernel is already in product form, but the increase in its dimension $k$ is essential.

## 5. Random regression.

Suppose that we observe an i.i.d. sample $(X_1, Y_1)$, $\ldots, (X_n, Y_n)$ from the distribution of a vector $(X, Y)$ described structurally as $Y = f(X) + \varepsilon$, for $(X, \varepsilon)$ a random vector with $E(\varepsilon \mid X) = 0$ and $\sigma^2(x) = E(\varepsilon^2 \mid X = x)$ admitting a bounded version. The distribution $P_X$ of $X$ is known and $\theta_1, \theta_2, \ldots$ are the Fourier coefficients of the regression function $f$ relative to a given orthonormal basis $e_1, e_2, \ldots$ of $L_2(P_X)$. We assume that the set of regression functions is uniformly bounded.

Given an initial estimator $\hat{\theta}^{(n)}$ our estimator for $\|\theta - \hat{\theta}^{(n)}\|^2$ is given by

$$R_{k,n}(\hat{\theta}^{(n)}) = \frac{1}{n(n-1)} \sum_{r \neq s=1}^{n} \sum_{i=1}^{k} (Y_r e_i(X_r) - \hat{\theta}_i^{(n)})(Y_s e_i(X_s) - \hat{\theta}_i^{(n)}).$$

Here $k = k_n$ is chosen dependent on $\Theta$. We combine this with the variance estimator

$$\hat{\tau}_{k,n,\theta}^2 = \frac{2k(\|f\|_\infty^2 + \|\sigma^2\|_\infty)^2}{n(n-1)} + \frac{4\|f\|_\infty^2 + 4\|\sigma\|_\infty^2}{n} \sum_{i=1}^k (\theta_i - \hat{\theta}_i^{(n)})^2.$$

THEOREM 5.1. *For any $k, n$,*

$$\sup_{\theta \in \Theta} \mathrm{E}_\theta \left( \left( \frac{R_{k,n}(\hat{\theta}^{(n)}) - \sum_{i=1}^k (\theta_i - \hat{\theta}_i^{(n)})^2}{\hat{\tau}_{k,n,\theta}} \right)^2 \Big| \hat{\theta}^{(n)} \right) \leq 1.$$

PROOF. The proof is similar to the proof of Theorem 4.1. The variable $R_{k,n}(\hat{\theta}^{(n)})$ is again a $U$-statistic of order 2. It has mean $\sum_{i=1}^k (\theta_i - \hat{\theta}_i^{(n)})^2$ and Hoeffding decomposition [cf. (4.1), but replace $X_i$ by $(X_i, Y_i)$] with kernels of the form

$$P_1 h(x, y) = 2 \sum_{i=1}^k (\theta_i - \hat{\theta}_i^{(n)})(y e_i(x) - \theta_i),$$

$$P_{1,2} h(x_1, y_1, x_2, y_2) = \sum_{i=1}^k (y_1 e_i(x_1) - \theta_i)(y_2 e_i(x_2) - \theta_i)$$

$$= \sum_{i=1}^k y_1 y_2 e_i(x_1) e_i(x_2)$$

$$- \sum_{i=1}^k \theta_i(y_1 e_i(x_1) + y_2 e_i(x_2)) + \sum_{i=1}^k \theta_i^2.$$

By the orthonormality of the functions $e_i$ and arguments as in the proof of Theorem 4.1,

$$\mathrm{var}\, P_1 h(X, Y) \leq 4 \|\mathrm{E}(Y^2|X)\|_\infty \sum_{i=1}^k (\theta_i - \hat{\theta}_i^{(n)})^2,$$

$$\mathrm{var}\, P_{1,2} h(X_1, Y_1, X_2, Y_2) \leq \|\mathrm{E}(Y^2|X)\|_\infty^2 k.$$

From $Y = f(X) + \varepsilon$ and $\mathrm{E}(\varepsilon|X) = 0$ it follows that $\mathrm{E}(Y^2|X) = f^2(X) + \mathrm{E}(\varepsilon^2|X) \leq \|f\|_\infty^2 + \|\sigma^2\|_\infty$. Combining the preceding bounds we obtain the theorem. □

The bound given by the preceding theorem is of the same form as the bounds given in the preceding sections, but with $\|f\|_\infty^2 + \|\sigma^2\|_\infty$ playing the role of $\sigma^2$ in Section 3. Again (2.8) is justified with $\hat{\tau}_{k,n,\theta}^2$ of the order as in (2.5). Proposition 2.1 gives the same corollaries for confidence regions.

**6. Lower bounds.** In this section we relate the minimum diameter of a confidence region to the minimax rates for testing and estimation. Consider a sequence of statistical experiments $(P_\theta^{(n)} : \theta \in \Theta)$ indexed by a parameter $\theta \in \Theta$ in a metric space $(\Theta, d)$ and a submodel indexed by a subset $\Theta_1 \subset \Theta$. We are interested in the maximal diameter over $\Theta_1$ of confidence regions that are honest over the whole model $\Theta$.

We shall silently understand that appropriate measurability assumptions regarding the confidence regions are satisfied.

Given $0 < \alpha < \beta < 1$, let $\varepsilon_n$ be a sequence of positive numbers such that there exists no sequence of tests $\phi_n$ satisfying the two requirements, for some given subsets $\Theta_{n,1} \subset \Theta_1$,

$$(6.1) \qquad \limsup_{n \to \infty} \sup_{\theta \in \Theta \,:\, d(\theta, \Theta_{n,1}) > \varepsilon_n} P_\theta^{(n)} \phi_n < \alpha,$$

$$(6.2) \qquad \limsup_{n \to \infty} \sup_{\theta \in \Theta_{n,1}} P_\theta^{(n)} (1 - \phi_n) < \beta.$$

This can only be satisfied if $\alpha + \beta \leq 1$, because otherwise the trivial test $\phi_n \equiv \alpha'$ for some $\alpha'$ with $\alpha' < \alpha$ and $1 - \alpha' < \beta$ satisfies (6.1)–(6.2). For $\beta \leq 1 - \alpha < 1$, the condition is satisfied for $\varepsilon_n$ equal to what Ingster [23] calls a rate of "not asymptotic indistinguishability of the hypotheses." The following lemma shows that the diameter over $\Theta_1$ of an honest confidence set is at least of the order $\varepsilon_n$.

LEMMA 6.1. *For given $0 < \alpha < \beta < 1$ and subsets $\Theta_{n,1} \subset \Theta_1$, if there exists no sequence of tests $\phi_n$ satisfying (6.1)–(6.2), then for any sequence of confidence sets $\hat{C}_n$ satisfying (1.1),*

$$\limsup_{n \to \infty} \sup_{\theta \in \Theta_1} P_\theta^{(n)} (\mathrm{diam}(\hat{C}_n) \geq \varepsilon_n) > \beta - \alpha.$$

PROOF. Let $\Theta_{n,0} = \{\theta \in \Theta : d(\theta, \Theta_{n,1}) > \varepsilon_n\}$. Given a sequence of confidence sets $\hat{C}_n$ satisfying (1.1) define tests by $\phi_n = \mathbb{1}_{d(\hat{C}_n, \Theta_{n,0}) > 0}$.

If $\theta \in \Theta_{n,0}$ and $d(\hat{C}_n, \Theta_{n,0}) > 0$, then $\theta \notin \hat{C}_n$. Therefore, from (1.1) it is immediate that these tests satisfy (6.1).

If $\theta \in \Theta_{n,1}$, $d(\hat{C}_n, \Theta_{n,0}) = 0$ and $\theta \in \hat{C}_n$, then $\mathrm{diam}(\hat{C}_n) \geq \varepsilon_n$. [Indeed, for every $\delta > 0$ there exist points $c \in \hat{C}_n$ and $\theta_n \in \Theta_{n,0}$ with $d(c, \theta_n) < \delta$. By the definition of $\Theta_{n,0}$ we have $d(\theta_n, \Theta_{n,1}) > \varepsilon_n$ and hence $d(\theta_n, \theta) > \varepsilon_n$. By the triangle inequality $d(c, \theta) > \varepsilon_n - \delta$.] It follows that, for every $\theta \in \Theta_{n,1}$,

$$P_\theta^{(n)} (1 - \phi_n) = P_\theta^{(n)} (d(\hat{C}_n, \Theta_{n,0}) = 0)$$
$$\leq P_\theta^{(n)} (\mathrm{diam}(\hat{C}_n) \geq \varepsilon_n) + P_\theta^{(n)} (\theta \notin \hat{C}_n).$$

By (1.1) the second term on the right-hand side is strictly asymptotically smaller than $\alpha$, uniformly in $\theta \in \Theta$. If the first term on the right-hand side were asymptotically smaller than $\beta - \alpha$, uniformly in $\theta \in \Theta_1$, thus contradicting the assertion of the lemma, then the left-hand side would be asymptotically strictly less than $\beta$, so that the tests would also satisfy (6.2). □

To obtain a lower bound for $\sup_{\theta \in \Theta_1} P_\theta^{(n)}(\operatorname{diam}(\hat{C}_n) > \varepsilon_n)$ we can apply the preceding lemma with $\Theta_{n,1} = \Theta_1$, but also with every subset of $\Theta_1$. In particular, we may apply the lemma with a one-point set $\Theta_{n,1} = \{\theta_1\}$, for any $\theta_1 \in \Theta_1$. For regularity models $\Theta$, Ingster [23] characterizes the minimax rate for exactly these one-point problems. He shows that there exists a rate $\varepsilon_n^*$ such that the sum of the error probabilities (6.1)–(6.2) goes to zero if $\varepsilon_n/\varepsilon_n^* \to \infty$ and goes to 1 if $\varepsilon_n/\varepsilon_n^* \to 0$. Thus the condition of the lemma is satisfied for any $0 < \alpha < \beta < 1$ with $\alpha + \beta \leq 1$ and $\varepsilon_n$ with $\varepsilon_n/\varepsilon_n^* \to 0$. The lemma then says that the weak limit points in $[0, \infty]$ of the distribution of $\operatorname{diam}(\hat{C}_n)/\varepsilon_n^*$ have a component of size at least $\beta - \alpha$ concentrated on $(0, \infty]$. In other words, the order of the diameter is at least $\varepsilon_n^*$.

The relationship between the diameter of confidence regions and the minimax rate for estimation is less perfect, due to the fact that the risk for estimation concerns the complete distribution of an estimator, whereas a confidence region at level $1 - \alpha$ leaves a mass of size $\alpha$ completely undiscussed.

A key result is as follows. Let $\beta \geq 0$ be given, and let $\varepsilon_n$ be a sequence of positive numbers such that for every estimator sequence $T_n$

$$(6.3) \qquad \liminf_{n \to \infty} \sup_{\theta \in \Theta_1} P_\theta^{(n)}(d(T_n, \theta) \geq \varepsilon_n) > \beta.$$

LEMMA 6.2. *For given $0 < \alpha < \beta < 1$, if (6.3) holds for every estimator sequence $T_n$, then for any sequence of confidence sets $\hat{C}_n$ satisfying (1.1),*

$$\liminf_{n \to \infty} \sup_{\theta \in \Theta_1} P_\theta^{(n)}(\operatorname{diam}(\hat{C}_n) \geq \varepsilon_n) > \beta - \alpha.$$

PROOF. Given a sequence of confidence sets $\hat{C}_n$, define for each $n$ an estimator $T_n$ to be an arbitrary point in $\hat{C}_n$. Then, for any $\theta \in \Theta_1$,

$$P_\theta^{(n)}(d(T_n, \theta) \geq \varepsilon_n) \leq P_\theta^{(n)}(\operatorname{diam}(\hat{C}_n) \geq \varepsilon_n) + P_\theta^{(n)}(\theta \notin \hat{C}_n).$$

By (1.1) the second term on the right-hand side is asymptotically smaller than $\alpha$, uniformly in $\theta \in \Theta$. By assumption the lim inf of the supremum of the left-hand side over $\theta \in \Theta_1$ is bounded below by $\beta$. □

If we choose $\varepsilon_n$ faster than the minimax rate, then typically (6.3) holds for some $\beta > 0$. In particular this is true if the minimax rate $\varepsilon_n^*$ has the property that for a "best" estimator sequence $T_n$ the sequence $d(T_n, \theta)/\varepsilon_n^*$ has all its limit points on $(0, \infty]$. In that case $d(T_n, \theta)/\varepsilon_n \to \infty$, and the right-hand side of (6.3) is 1, for any sequence $\varepsilon_n$ with $\varepsilon_n/\varepsilon_n^* \to 0$. We may then apply the lemma with any $\beta < 1$. More generally, this argument works if the weak limit points of the sequence $d(T_n, \theta)/\varepsilon_n^*$ in $[0, \infty]$ possess a point mass of at most $\beta$ at 0.

DEPARTMENT OF EPIDEMIOLOGY
HARVARD SCHOOL OF PUBLIC HEALTH
677 HUNTINGTON AVENUE
BOSTON, MASSACHUSETTS 02115
USA
E-MAIL: robins@hsph.harvard.edu

DEPARTMENT OF MATHEMATICS
VRIJE UNIVERSITEIT AMSTERDAM
DE BOELELAAN 1081A
1081 HV AMSTERDAM
THE NETHERLANDS
E-MAIL: aad@cs.vu.nl